\documentclass[twoside,leqno]{article}
\usepackage{amsmath}

\catcode`\@=11
\def\keywords#1{\def\@keywords{#1}}
\let\@keywords=\@empty
\def\subjclass#1{\def\@subjclass{#1}}
\let\@subjclass=\@empty
\newcommand{\keywordsname}{Key words and phrases}
\newcommand{\subjclassname}{\textup{2000} Mathematics Subject
     Classification}
\def\@addpunct#1{\ifnum\spacefactor>\@m \else#1\fi}
\def\@setsubjclass{%
  {\itshape\subjclassname.}\enspace\@subjclass\@addpunct.}
\def\@setkeywords{%
  {\itshape \keywordsname.}\enspace \@keywords\@addpunct.}
\renewcommand\thesection       {\arabic{section}}

\def\initrhead{}%
\def\initrhead{{\it JAMS}}%
\newcounter{footpage}%
\def\sikibetuh{}%
\mark{{}{}}%
\def\ps@myheadings{%
    \let\@oddfoot\@empty\let\@evenfoot\@empty
    \def\@evenhead{\small\sikibetuh\thepage\hfill\leftmark\hfill}%
    \def\@oddhead{\small\hfill\rightmark\hfill\sikibetuh\thepage}%
    \let\@mkboth\@gobbletwo
    \let\sectionmark\@gobble
    \let\subsectionmark\@gobble
    }
\def\ps@jamsinit{%
    \def\@evenhead{\small\sikibetuh\thepage\hfill{\footnotesize{\initrhead}}}%
    \def\@oddhead{\small{\footnotesize{\initrhead}}\hfill\sikibetuh\thepage}%
    \let\@mkboth\@gobbletwo
    \let\sectionmark\@gobble
    \let\subsectionmark\@gobble
    }
\def\section{\@startsection{section}{1}%
  \z@{.7\baselineskip\@plus\baselineskip}{-.5em}%
  {\normalfont\bfseries}}
\def\subsection{\@startsection{subsection}{2}%
  \z@{.5\baselineskip\@plus.7\baselineskip}{-.5em}%
  {\normalfont\bfseries}}
\def\subsubsection{\@startsection{subsubsection}{3}%
  \z@{.5\baselineskip\@plus.7\baselineskip}{-.5em}%
  {\normalfont\itshape}}
\def\paragraph{\@startsection{paragraph}{4}%
  \z@\z@{-\fontdimen2\font}%
  \normalfont}
\def\subparagraph{\@startsection{subparagraph}{5}%
  \z@\z@{-\fontdimen2\font}%
  \normalfont}
\def\appendix{\par\c@section\z@ \c@subsection\z@
   \let\sectionname\appendixname
   \def\thesection{\@Alph\c@section}}
\def\appendixname{Appendix}
\renewenvironment{abstract}{%
\noindent
  \small
\noindent
  \quote
  \parindent=1.5em
  {\scshape \abstractname .\hspace{\z@}}%
  }
{  \endquote
}
\renewcommand{\maketitle}{\par
  \begingroup
    \renewcommand{\thefootnote}{\fnsymbol{footnote}}%
    \def\@makefnmark{\hbox to\z@{$\m@th^{\@thefnmark}$\hss}}%
    \long\def\@makefntext##1{\parindent 1em\noindent
            \hbox to1.8em{\hss$\m@th^{\@thefnmark}$}##1}%
      \newpage
      \global\@topnum\z@   %
      \@maketitle
    \@thanks
  \endgroup
  \thispagestyle{jamsinit}
  \setcounter{footnote}{0}%
  \let\thanks\relax
  \let\maketitle\relax\let\@maketitle\relax
  \gdef\@thanks{}
  \gdef\@author{}\gdef\@title{}}
\def\@maketitle{%
  \newpage
  \null
  \vskip 2em%
  \begin{center}%
    {\normalfont\bfseries \@title \par}%
    \vskip 1.5em%
    {\normalfont\scshape
      \lineskip .5em%
      \begin{tabular}[t]{c}%
        \@author
      \end{tabular}\par
     \def\thefootnote{}
     \ifx\@empty\@subjclass\else \footnotetext{\@setsubjclass}\fi
     \def\thefootnote{}
     \ifx\@empty\@keywords\else \footnotetext{\@setkeywords}\fi
}%
    \vskip 1em%
    {\normalfont \@date}%
  \end{center}%
  \par
  \vskip 0.8em
}
\renewenvironment{thebibliography}[1]
      {
  \@startsection{section}
        \@m\z@{9\p@\@plus12\p@}{6\p@}%
        {\centering\scshape}\refname
  \small\labelsep .5em\relax
      \list{\@biblabel{\@arabic\c@enumiv}}%
           {\settowidth\labelwidth{\@biblabel{#1}}%
            \leftmargin\labelwidth
            \advance\leftmargin\labelsep
            \@openbib@code
            \usecounter{enumiv}%
            \let\p@enumiv\@empty
            \renewcommand\theenumiv{\@arabic\c@enumiv}}%
      \sloppy
      \clubpenalty4000
      \@clubpenalty \clubpenalty
      \widowpenalty4000%
      \sfcode`\.\@m}
     {\def\@noitemerr
       {\@latex@warning{Empty `thebibliography' environment}}%
      \endlist}

\usepackage{graphicx}
\usepackage{amsfonts}
\usepackage{amssymb}
\newtheorem{theorem}{Theorem}
\newtheorem{corollary}[theorem]{Corollary}
\newtheorem{definition}[theorem]{Definition}

\newtheorem{proposition}[theorem]{Proposition}
\newtheorem{remark}[theorem]{Remark}
\begin{document}
\title{Lusin's Theorem and Bochner Integration}
\author{Peter A. Loeb\\Department of Mathematics, University of 
Illinois\\1409 West Green Street, 
Urbana, Illinois 61801, U.S.A.\\e-mail: loeb@math.uiuc.edu
\and Erik Talvila\thanks{The work of both authors here was supported in part by the NSERC and the University of Alberta. The authors thank David Ross for a helpful pointer to the literature in \cite{fremlin0} and \cite{kupkaprikry}.}\\Department of Mathematics and Statistics \\ University College of the Fraser Valley\\Abbotsford, BC, Canada V2S 7M8\\e-mail: Erik.Talvila@ucfv.ca}
\date{To appear in {\it  Scientiae Mathematicae Japonicae.}  
Accepted January 13, 2004.}

\maketitle

\begin{abstract}
It is shown that the approximating functions used to define the Bochner
integral can be formed using geometrically nice sets, such as balls, from a
differentiation basis. Moreover, every appropriate sum of this form will be
within a preassigned $\varepsilon$ of the integral, with the sum for the local
errors also less than $\varepsilon$. All of this follows from the ubiquity of
Lebesgue points, which is a consequence of Lusin's theorem, for which a simple
proof is included in the discussion.
\end{abstract}

\section{Introduction}

An attractive feature of the Henstock--Kurzweil integral for the real line is
that it can be defined in terms of Riemann sums using intervals. In
\cite{loebtalvila}, we showed that for more general spaces, Lebesgue points
and points of approximate continuity can be used to approximate a Lebesgue
integral $\int_{X}f~d\mu$ with sums employing disjoint, geometrically nice
sets $S_{i}$ covering all but a set of measure $0$ of $X$. We also showed that
any such sum will be within a preassigned $\varepsilon$ of the integral,
provided each term $f(x_{i})\mu(S_{i})$ of the sum has the property that
$S_{i}$ is contained in a ball about $x_{i}$ of radius at most $\delta(x_{i})$, 
where $\delta$ is a \textquotedblleft gauge\textquotedblright\ function
determined by $f$ and $\varepsilon$.

In this article, we show that a similar result holds for the Bochner integral,
where the domain of the integrand is a measure space 
$\left(  X,\mathcal{M},\mu\right)$ on which is defined a differentiation basis. 
Although somewhat more general settings are possible, we will assume that $X$ is a finite
dimensional normed vector space and the differentiation basis is obtained
using the Besicovitch or Morse covering theorem (see Section \ref{Covering}).
The geometrically nice sets we will use for the approximating functions will
be balls or the more general starlike sets described in Section~\ref{Covering}. 
We will assume that $\mu$ is a complete \textbf{Radon measure}; i.e., $\mu$
is a regular measure on $\mathcal{M}$, which includes the Borel sets, and
compact sets have finite measure. By a \textbf{measurable} set, we will always
mean a set in $\mathcal{M}$. We will let $\mathbb{N}$ denote the natural
numbers and $\mathbb{R}$ the real numbers.

In what follows, an integrand $f$ will take its values in a Banach space
$\left(  Y,\left\Vert \cdot\right\Vert \right)  $ and will be 
$\mu$\textbf{-measurable}, meaning, there is a sequence of simple functions
$f_{n}$ with $\lim_{n}\left\Vert f_{n}-f\right\Vert =0$ $\mu$-a.e. on $X$.
This will allow us to use Lusin's theorem. An elementary proof of that
theorem, simple even for real-valued functions on $\mathbb{R}$, is given in
the next section.

Our approximating sums will be constructed using a gauge function $\delta$
mapping $X$ into $(0,1]$. In the theory of Henstock--Kurzweil and McShane
integration, the appearance of a gauge function is somewhat mysterious. We
show in proving Theorem \ref{integration} how the properties of Lebesgue
points (discussed below in Section \ref{Covering}) can be used to determine an
appropriate gauge.

We note that the Bochner integral has been studied in terms of Riemann sums
over finite partitions of a metric space in \cite{Neerven} and over
generalized McShane partitions of a measure space in \cite{fremlin} and
\cite{dipiazza}. See also \cite{dipiazzamusial}. None of these papers,
however, works with approximations using geometrically nice sets.

\section{Lusin's Theorem}

\label{LusinSec}

Recall that a function $f$ from $X$ into a topological space $(Y,\mathcal{T})$
is \textbf{Borel measurable} if the inverse image of each open set in $Y$ is
in $\mathcal{M}$. In our case, where $Y$ is a Banach space and $f$ is $\mu
$-measurable, it follows from Theorem III.6.10 of \cite{dunfordschwartz} that
$f$ is Borel measurable on $X$. We note that it follows from the same theorem
and a deep result of D. H. Fremlin (Theorem 2B in \cite{fremlin0} or Theorem
4.1 in the expository article \cite{kupkaprikry} by J. Kupka and K. Prikry)
that if $Y$ is a metric space and  $f:X\rightarrow Y$ is Borel measurable,
then $f$ is $\mu$-measurable.

In any case, Lusin's theorem holds for the restriction of $f$ to a set
$\Omega\subseteq X$ with $\mu\left(  \Omega\right)  <+\infty$, and is easily
extended to all of $X$ using the $\sigma$-finiteness of $\mu$. Here is the
general statement of the theorem for the case of a finite measure space, with
an elementary proof that is appropriate even for the simplest setting.

\begin{theorem}
[Lusin]\label{Lusin}Let $Y$ be a topological space with a countable base
$\left\langle V_{n}\right\rangle $ for the topology, and let $f$ be a Borel
measurable function from a Radon measure space of finite measure
$(\Omega,\mathcal{A},\mu)$ into $Y$. Given $\varepsilon>0$, there is a compact
set $K$ with $\mu(\Omega\setminus K)<\varepsilon$ such that $f$ restricted to
$K$ is continuous.
\end{theorem}

\noindent \textbf{Proof:} Fix compact sets $K_{n}\subseteq f^{-1}[V_{n}]$ and $K_{n}^{\prime}%
\subseteq\Omega\setminus f^{-1}[V_{n}]$ for each $n$ so that $\mu
(\Omega\setminus K)<\varepsilon$ when $K:=\bigcap\nolimits_{n}(K_{n}\cup
K_{n}^{^{\prime}})$. Given $x\in K$ and an $n$ with $f(x)\in V_{n}$,\ \ $x\in
O:=\Omega\setminus K_{n}^{\prime}$ and $f\left[  O\cap K\right]  \subseteq
V_{n}$. $\square$

\begin{remark}
\begin{rm}
In Oxtoby's text \cite{oxtoby} a similar principle is used in a more complex
proof to show that a measurable $f:\mathbb{R}\rightarrow\mathbb{R}$ is
continuous when restricted to a large measurable subset of $\mathbb{R}$.
\end{rm}
\end{remark}

\section{Covering Theorems}

\label{Covering}

Our integration result is based on calculations using a covering theorem. Here
we use either Besicovitch's theorem \cite{besicovitch} (also
see\ \cite{furediloeb}) for a covering by balls, or the theorem of Morse
\cite{morse} involving more general sets. In \cite{loebtalvila}, we have
established strengthened versions of these theorems that hold for our finite
dimensional normed vector space $X$. Moreover, these covering theorems are
also valid for a space that is locally isometric to $X$ since one only needs
to bound the cardinality of a finite collection of sets with small diameter
forming a \textquotedblleft$\tau$-satellite configuration\textquotedblright 
as defined in \cite{furediloeb} and \cite{loebtalvila}. We leave it to the
interested reader to consider our results in this more general setting, as
well as settings using Vitali's covering theorem.

We will denote the closed (and compact) ball in $X$ with center $a$ and radius
$r>0$ by $B(a,r):=\{x\in X:\Vert x-a\Vert\leq r\}$. The sets used by Morse
involve a parameter $\lambda\geq1$; they are closed balls when $\lambda=1$.

Given $\lambda\geq1$ and $a\in X$, we say that a set $S(a)\subseteq X$ is a
\textbf{Morse set} or $\lambda$\textbf{-Morse set} associated with $a$ and
$\lambda$ if it satisfies two conditions. First, $S(a)$ must be 
$\lambda$\textbf{-regular}. This means that there is an $r>0$ such that
$B(a,r)\subseteq S(a)\subseteq B(a,\lambda r)$. Second, $S(a)$ must be
\textbf{starlike} with respect to $B(a,r)$. This means that for each 
$y \in B(a,r)$ and each $x \in S(a)$, the line segment 
$\alpha y+(1-\alpha)x$,\ \ $0\leq\alpha\leq1$, is contained in $S(a)$. 
We call $a$ the \textbf{tag} for $S(a)$. Note that the closure 
$\operatorname*{cl}\left(  S(a)\right)$ of a
$\lambda$-Morse set $S(a)$ is again a $\lambda$-Morse set. We say that a
$\lambda$-Morse set $S(a)$ is $\delta$\textbf{-fine} with respect to a gauge
function $\delta$ defined at $a$ if $S(a)\subseteq B(a,\delta(a))$.

Suppose we are given $\lambda\geq1$, a Radon measure $\mu$ on $X$, and an open
subset $\Omega$ of $X$. We will call a collection $\mathcal{S}$ of 
$\lambda$-Morse sets a \textbf{fine},\textbf{ measurable}, 
$\lambda$\textbf{-Morse cover }of $\Omega$ provided each set in 
$\mathcal{S}$ is a measurable subset
of $\Omega$ and each $a\in\Omega$ is the tag of sets in $\mathcal{S}$ with
arbitrarily small diameters. A sequence $\left\langle S_{i}\right\rangle $
from such an $\mathcal{S}$ is said to be $\mu$\textbf{-exhausting of} $\Omega$
if it is a finite or countably infinite sequence that is pairwise disjoint,
and covers all but a set of $\mu$-measure $0$ of $\Omega$. A collection
$\mathcal{S}$ is called a $\mu$\textbf{-a.e.\hspace{0.03in}}$\lambda
$\textbf{-Morse} \textbf{cover} of $\Omega$ if it is a fine, measurable,
$\lambda$-Morse cover of $\Omega$ and for each nonempty open subset
$U\subseteq\Omega$ the collection 
$\mathcal{S}_{U}=\left\{  S\in \mathcal{S}:S\subseteq U\right\}  $ 
has the property that for any gauge
function $\delta:U\rightarrow(0,1]$, there is a sequence of $\delta$-fine sets
in $\mathcal{S}_{U}$ that is $\mu$-exhausting of $U$. In \cite{loebtalvila} we
have established the following result for $\mu$ and $\Omega$.

\begin{proposition}
\label{prop1} A fine, measurable, $\lambda$-Morse cover $\mathcal{S}$ \ of
$\Omega$ is a $\mu$-a.e. $\lambda$-Morse\hspace{0.02in}cover of $\Omega$ if it
consists of closed sets or if for each set 
$S\in\mathcal{S}$, $\mu(\Omega \cap(\operatorname*{cl}(S)\setminus S))=0$.
\end{proposition}

We have also shown that the conditions in Proposition~\ref{prop1} are
fulfilled by any measurable Morse cover $\mathcal{S}$ that is \textbf{scaled}.
This means that for each $S(a)\in\mathcal{S}$ and each $p\in(0,1]$, the set
$S^{(p)}(a)$ is also in $\mathcal{S}$ where $S^{(p)}(a)=\{a+px:a+x\in S(a)\}$.

\section{Approximate Continuity and Lebesgue Points}

\label{LebPt}

For this section, we fix a $\mu$-a.e. $\lambda$-Morse cover $\mathcal{S}$ of
an open set $\Omega\subseteq X$. We work with a $\mu$-measurable
$f:\Omega\rightarrow Y$. The following notions depend on the choice of
$\mathcal{S}$.

\begin{definition}
A point $a\in\Omega$ is a \textbf{point of approximate continuity} for $f$ if
for all positive $\varepsilon$ and $\eta$ there is an $R>0$ such that if
$S(a)$ is a set in $\mathcal{S}$ with tag $a$ and $S(a)\subseteq B(a,R)$, then
for $E(a,\eta):=\{x\in S(a):\left\Vert f(a)-f(x)\right\Vert >\eta\}$ we have
$\mu(E(a,\eta))\leq\varepsilon\,\mu(S(a))$. A point $a\in\Omega$ is a
\textbf{Lebesgue point} of $f$ with respect to $f(a)$ if for any
$\varepsilon>0$ there is an $R>0$ such that if $S(a)$ is a set in
$\mathcal{S}$ with tag $a$ and $S(a)\subseteq B(a,R)$, then
\[
\int\nolimits_{S(a)}\left\Vert ~\!\!f(x)-f(a)\right\Vert \,\mu(dx)\leq
\varepsilon\,\mu(S(a)).
\]

\end{definition}

It is easy to see that if $a\in\Omega$ is a Lebesgue point of $f$ with respect
to $f(a)$, then $a$ is a point of approximate continuity for $f$. It follows
from the fact that $\mathcal{S}$ is a differentiation basis that if $g$ is a
$\mu$-integrable, nonnegative, real-valued function on $\Omega$, then $\mu
$-almost all points of $\Omega$ are Lebesgue points. (See, for example,
\cite{bliedtnerloeb}.) Moreover, if $A$ is a measurable subset of $\Omega$,
then almost all points of $A$ are \textbf{points of density}, that is, points
of approximate continuity with respect to the characteristic function
$\chi_{A}$ of $A$. Therefore, we have the following consequence of Lusin's theorem.

\begin{proposition}
If $f:\Omega\rightarrow Y$ is $\mu$-measurable, then $\mu$-almost all points
of $\Omega$ are points of approximate continuity for $f$.
\end{proposition}

\noindent \textbf{Proof:} If $\mu\left(  \Omega\right)  <+\infty$, then by Lusin's theorem 
(Theorem \ref{Lusin}), there is an increasing sequence of compact sets 
$K_{n}\subseteq\Omega$ such that for each $n$, $f|K_{n}$ is continuous and
$\mu\left(  \Omega\backslash\cup_{n}K_{n}\right)  =0$. For this case, the
result follows from the fact that for each $n$, $\mu$-almost every point of
$K_{n}$ is a point of density of $K_{n}$. The general case follows since
$\Omega$ has $\sigma$-finite measure. $\square$

\bigskip

We also have the following relationship between points of approximate
continuity and Lebesgue points.

\begin{proposition}
\label{lebisleb}If $f:\Omega\rightarrow Y$ is $\mu$-measurable and $a\in
\Omega$ is a point of approximate continuity for $f$ and also a Lebesgue point
for $\left\Vert f\right\Vert $ with respect to $\left\Vert f(a)\right\Vert $,
then $a$ is a Lebesgue point for $f$ with respect to $f(a)$.
\end{proposition}

\noindent \textbf{Proof:} Fix $\varepsilon>0$. Let $c=\left\Vert f(a)\right\Vert $ 
and choose $R>0$ so that if $S$ is a set in $\mathcal{S}$ with tag $a$ and 
$S\subseteq B(a,R)$, then
\[
\int_{S}\left\vert ~\left\Vert f(x)\right\Vert -c~\right\vert \ \mu
(dx)<\varepsilon\cdot\mu\left(  S\right)  ,
\]
and for $E:=\{x\in S:\left\Vert f(x)-f(a)\right\Vert >\varepsilon\}$, we have
$\mu(E)\leq\frac{\varepsilon}{2c+1}\cdot\mu\left(  S\right)  $. Now
\[
\int_{E\cap\left\{  \left\Vert f\right\Vert >c\right\}  }\left(  \left\Vert
f(x)\right\Vert -c\right)  \ \mu(dx)<\varepsilon\cdot\mu\left(  S\right)  ,
\]
whence
\[
\int_{E\cap\left\{  \left\Vert f\right\Vert >c\right\}  }\left\Vert
f(x)\right\Vert \ \mu(dx)\leq2\varepsilon\cdot\mu\left(  S\right)  .
\]
Therefore,
\begin{align*}
&  \int_{S}\left\Vert f(x)-f(a)\right\Vert \ \mu(dx)\\
&  \leq\int_{S\backslash E}\varepsilon\ \mu(dx)+\int_{E}2c\ \mu(dx)+\int
_{E\cap\left\{  \left\Vert f\right\Vert >c\right\}  }\left\Vert
f(x)\right\Vert \ \mu(dx)\leq4\varepsilon\cdot\mu\left(  S\right) .
\end{align*}
Since the choice of $\varepsilon>0$ is arbitrary, the result follows.
$\square$

\bigskip

\begin{remark}
\begin{rm}
Proposition \ref{lebisleb} also holds for the $\mu$-null set consisting of
points $a\in\Omega$ that are Lebesgue points of $\left\Vert f\right\Vert $
with respect to values different from $\left\Vert f(a)\right\Vert $.%
\end{rm}
\end{remark}

\section{Integration}

\label{Integration}

Recall that if $\mu\left(  \Omega\right)  <+\infty$, then a $\mu$-measurable
function $f:\Omega\rightarrow Y$ is \textbf{Bochner integrable} (with respect
to $\mu$) if there is a sequence of simple functions 
$f_{n}:$ $\Omega \rightarrow Y$ with 
$\lim_{n}\int_{\Omega}\left\Vert f_{n}-f\right\Vert d\mu=0$. 
In this case, the sequence of integrals $\int_{\Omega}f_{n}~d\mu$ is
Cauchy in $Y$, and the limit is the Bochner integral of $f$. Also, $f$ \ is
Bochner integrable if and only if the norm $\left\Vert f\right\Vert $ is
integrable. (See, for example, Chapter II of \cite{diesteluhl}.) We will
consider the analogous approximations and integral for the case that
$\left\Vert f\right\Vert $ is integrable on our $\sigma$-finite measure space
$\left(  X,\mathcal{M},\mu\right)  $, with the understanding that a simple
function is defined on all of $X$ but must vanish off of a set of finite
measure. In Theorem 12 of \cite{loebtalvila}, we have established a necessary
and sufficient condition for the integrability of a real-valued function, such
as $\left\Vert f\right\Vert $, in terms of approximations by Riemann sums
using sets from a $\mu$-a.e. Morse cover. Now, assuming that $\left\Vert
f\right\Vert $ is $\mu$-integrable on $X$, we show that the Bochner integral
of $f$ can also be approximated by the integral of any appropriate simple
function using Morse sets.

\begin{theorem}
\label{integration}Assume $\left\Vert f\right\Vert $ is $\mu$-integrable on
$X$, and fix $\varepsilon>0$. Given $\lambda\geq1$, and any $\mu$-a.e.,
$\lambda$-Morse cover $\mathcal{S}$ \ of $\ X$, there is a gauge function
$\delta\!:\!X\rightarrow(0,1]$ such that for any finite or countably infinite
sequence $\left\langle S_{i}(x_{i})\right\rangle $ of $\delta$-fine sets from
$\mathcal{S}$ that is $\mu$-exhausting of $X$, the function 
$\sum_{i} f(x_{i})\chi_{S_{i}}$ approximates $f$ in the sense that
\[
\int_{X}\left\Vert f(y)-\sum\nolimits_{i}\,f(x_{i})\chi_{S_{i}}
\,(y)\right\Vert ~\mu\left(  dy\right)  <\varepsilon
\]
and the absolute sum of local errors is small; that is,
\begin{equation}
\sum\nolimits_{i}\left\Vert \int_{S_{i}}f(y)~\mu(dy)-f(x_{i})\mu
(S_{i})\right\Vert <\varepsilon. \label{localerrorsum}
\end{equation}
The same approximation holds for the restriction of $f$ \ to a large open ball
$\Omega$ about the origin with 
$\int_{X\diagdown\Omega}\left\Vert f\right\Vert <\varepsilon$, 
but then there is an $m\in\mathbb{N}$ such that for any $n\geq m$,
\begin{align*}
\lefteqn{\left\Vert \int_{X}f(y)~\mu\left(  dy\right)  -\sum\nolimits_{i=1}
^{n}\,f(x_{i})~\mu\left(  S_{i}\right)  \right\Vert }\\
&  \leq\int_{X}\left\Vert f(y)-\sum\nolimits_{i=1}^{n}\,f(x_{i})\chi_{S_{i}
}\,(y)\right\Vert ~\mu\left(  dy\right)  <3\varepsilon.
\end{align*}

\end{theorem}

\noindent \textbf{Proof:} Fix $\gamma>0$ so that for each $E\subseteq X$ with 
$\mu(E)<\gamma $,\ $\int_{E}\left\Vert f\right\Vert ~d\mu<\varepsilon/4$. 
Let $L$ be the set of points $x\in X$ that are Lebesgue points of $f$ with 
respect to $f(x)$, and let $A:=X\setminus L$. Since $\mu\left(  A\right)  =0$, 
there is an open set $G$ containing $A$ with $\mu(G)<\gamma$. 
For each $n\in\mathbb{N}$, we set
$A_{n}=\{x\in A:n-1\leq\left\Vert f(x)\right\Vert <n\}$. The sets $A_{n}$ are
disjoint and $\mu$-null with union $A$. For each $n\in\mathbb{N}$, fix an open
set $G_{n}$ with $G\supseteq$ $G_{n}\supseteq A_{n}$ and 
$\mu(G_{n})<\varepsilon/\left(  n\cdot2^{n+2}\right)  $. For each $x\in A_{n}$, we
choose $\delta(x)<1\wedge\sup\{s:B(x,s)\subseteq G_{n}\}$. Then for any finite
or countably infinite disjoint sequence of $\ \delta$-fine sets $S_{i}$ from
$\mathcal{S}$ with tags $\ x_{i}$ \ in $\ A$, we have
\[
\sum\limits_{i}\left\Vert f(x_{i})\right\Vert \,\mu(S_{i})\leq\sum
\limits_{n=1}^{\infty}\left(  n\sum\limits_{x_{i}\in A_{n}}\!\!\mu
(S_{i})\right)  \leq\sum\limits_{n=1}^{\infty}\varepsilon\,2^{-n-2}
=\varepsilon/4,
\]
and the integral of $\left\Vert f\right\Vert $ over $\cup_{i}S_{i}$ is at most
$\varepsilon/4$. Set $B(0,-1)=B(0,0)=\varnothing$. If for $n\in\mathbb{N}$,
$x\in L\cap\left(  B(0,n)\diagdown B(0,n-1)\right)  $, then $B(x,1)\subseteq
E_{n}:=B\left(  0,n+1\right)  \diagdown B\left(  0,n-2\right)  $, and we
choose $\delta(x)$ with $0<\delta(x)\leq1$ so that every $\delta$-fine 
$S\in\mathcal{S}$ with tag $x$ satisfies the inequality
\[
\int\nolimits_{S}\!\!\left\Vert f(y)-f(x)\right\Vert \,\mu(dy)<\frac
{\varepsilon\cdot2^{-n-2}}{\left[  1+\mu\left(  E_{n}\right)  \right]  }
\cdot\,\mu(S).
\]
With this choice of the gauge $\delta:X\rightarrow(0,1]$, we let $\left\langle
S_{i}\right\rangle $ be any finite or countably infinite sequence of $\delta
$-fine sets from $\mathcal{S}$ that is $\mu$-exhausting of $X$. Let $I_{L}$ be
the set of those indices $i$ for which $x_{i}\in L$, and let $I_{A}$ be the
set of those indices $i$ with $x_{i}\in A=X\setminus L$. For each
$n\in\mathbb{N}$, let $I_{L}^{n}$ be those indices in $I_{L}$ with $x_{i}\in
B(0,n)\backslash B(0,n-1)$. Set $S_{A}:=\cup_{i\in I_{A}}S_{i}(x_{i})$ and
$S_{L}:=\cup_{i\in I_{L}}S_{i}(x_{i})$. Now by the above calculation,
$\sum\limits_{i\in I_{A}}\left\Vert f(x_{i})\right\Vert \,\mu(S_{i}
)\leq\varepsilon/4$, and since $\mu\left(  X\backslash(S_{L}\cup
S_{A})\right)  =0$, we have
\[
\int_{X\backslash(S_{L}\cup S_{A})}\left\Vert f\right\Vert ~d\mu+\int_{S_{A}%
}\left\Vert f\right\Vert ~d\mu+\sum\limits_{i\in I_{A}}\left\Vert
f(x_{i})\right\Vert \,\mu(S_{i})\leq\frac{\varepsilon}{2}.
\]
Moreover,
\begin{align*}
&  \int_{X}\left\Vert f(y)-\sum\nolimits_{i}\,f(x_{i})\chi_{S_{i}
}\,(y)\right\Vert ~\mu\left(  dy\right) \\
&  \leq\int_{S_{L}}\left\Vert f(y)-\sum\nolimits_{i\in I_{L}}\,f(x_{i}
)\chi_{S_{i}}\,(y)\right\Vert ~\mu(dy)+\frac{\varepsilon}{2}\\
&  \leq\sum_{n=1}^{\infty}\sum_{i\in I_{L}^{n}}\int_{S_{i}}\!\!\left\Vert
f(y)-f(x_{i})\right\Vert ~\mu(dy)+\frac{\varepsilon}{2}\\
&  <\sum_{n=1}^{\infty}\sum_{i\in I_{L}^{n}}\frac{\varepsilon\cdot2^{-n-2}
}{\left[  1+\mu\left(  E_{n}\right)  \right]  }\cdot\,\mu(S(x_{i}
))+\frac{\varepsilon}{2}\leq\varepsilon\cdot\sum_{n=1}^{\infty}2^{-n-2}
+\frac{\varepsilon}{2}<\varepsilon.
\end{align*}

To show that the sum of local errors is small, we note that for each $i$,
$\sum\nolimits_{j}\,f(x_{j})\chi_{S_{j}}(y)=f(x_{i})$ on $S_{i}$, whence
\begin{align*}
&  \sum\nolimits_{i}\left\Vert \int_{S_{i}}f(y)~\mu(dy)-f(x_{i})\mu
(S_{i})\right\Vert \leq\sum\nolimits_{i}\int_{S_{i}}\left\Vert f(y)-f(x_{i}
)\right\Vert ~\mu(dy)\\
&  =\sum\nolimits_{i}\int_{S_{i}}\left\Vert f(y)-\sum\nolimits_{j}
\,f(x_{j})\chi_{S_{j}}(y)\right\Vert ~\mu(dy)\\
&  =\int_{X}\left\Vert f(y)-\sum\nolimits_{j}\,f(x_{j})\chi_{S_{j}
}(y)\right\Vert ~\mu(dy)<\varepsilon.
\end{align*}

For the last inequality, choose $m\in\mathbb{N}$ so that 
$\mu(\Omega \backslash\cup_{i\leq m}S_{i})<\gamma$, 
and the proof follows essentially as
before. $\ \square$

\begin{remark}
\begin{rm}
Note that our simple functions are not formed using a partition of $X$, as is
usual for the Henstock--Kurzweil and McShane integrals; they are, however,
defined on all of $X$. Also note that partitions are allowed in Theorem
\ref{integration}. That is, suppose a Cartesian coordinate system has been
imposed on $X$ and $I\subset X$ is a bounded open interval. In McShane
integration over $I$, one takes a partition 
$\left\langle S_{i}\right\rangle_{i=1}^{N}$ to be intervals in $I$ with 
$\cup_{i=1}^{N}S_{i}=I$. If these
intervals are mutually disjoint, satisfy the $\lambda$-regularity condition,
and are $\delta$-fine, then they are examples of the covering sets in Theorem
\ref{integration}.
\end{rm}
\end{remark}

\begin{definition}
A function $G:\mathcal{M}\rightarrow Y$ is called \textbf{countably additive} if
$G(\cup_{i}E_{i})=\lim_{n\rightarrow\infty}\sum_{i}^{n}G(E_{i})$ in norm for
every ordering of any infinite, pairwise disjoint sequence 
$\langle E_{i}\rangle$ from $\mathcal{M}$. It is also required that 
$G(\varnothing)=0$, whence $G$ is finitely additive.
\end{definition}

\begin{corollary}
\label{sumimpliesintegrability}Let $\mu$ be a Radon measure on 
$(X,\mathcal{M})$. Let $\mathcal{S}$ be a $\mu$-a.e., 
$\lambda$-Morse cover of $X$ such that
for each $S\in\mathcal{S}$, $\ \mu(S\diagdown S^{\circ})=0$, where $S^{\circ}$
is the interior of $S$. Suppose $G:\mathcal{M}\rightarrow Y$ is countably
additive and has the additional properties that $G(E)=0$ for each
$E\in\mathcal{M}$ with $\mu(E)=0$, and
\[
M:=\sup\left\{  \sum_{i}\left\Vert G(S_{i})\right\Vert :\left\langle
S_{i}\right\rangle \subset\mathcal{S}\text{,\ \ }\left\langle S_{i}
\right\rangle \text{\ }\mu\text{-exhausting of }X\right\}  <+\infty.
\]
Also suppose that there is a $\mu$-measurable $f:X\rightarrow Y$ such that for
any $\varepsilon>0$ there is a gauge $\delta:X\rightarrow(0,1]$ \ with the
property that for any $\delta$-fine sequence 
$\left\langle S_{i}(x_{i})\right\rangle $ in $\mathcal{S}$ 
that is $\mu$-exhausting of $X$ we have
\begin{equation}
\sum\nolimits_{i}\left\Vert f(x_{i})\mu(S_{i})-G(S_{i})\right\Vert
<\varepsilon. \label{riemannsum}
\end{equation}
Then $f$ is $\mu$-integrable and $\int_{X}f\,d\mu=G(X)$.
\end{corollary}

\noindent \textbf{Proof:} Fix $\varepsilon>0$. Also, fix 
$\langle T_{j}\rangle\subset\mathcal{S}$
so that $\langle T_{j}\rangle$ is a sequence that is $\mu$-exhausting of $X$
and $\sum_{j}\Vert G(T_{j})\Vert>M-\varepsilon/2$. Let 
$\Omega=\cup_{j} T_{j}^{\circ}$. We choose a gauge $\delta:X\rightarrow(0,1]$ 
so that\ for any $\delta$-fine sequence $\left\langle S_{i}\right\rangle $ in 
$\mathcal{S}$ that is $\mu$-exhausting of $X$, Equation \ref{riemannsum} holds with
$\varepsilon$ replaced by $\varepsilon/2$ \ and if $x\in T_{j}^{\circ}$ for
some $j$, then $B(x,\delta(x))\subseteq T_{j}^{\circ}$. The collection
$\mathcal{S}_{\Omega}$ of $\delta$-fine sets in $\mathcal{S}$ with tags in
$\Omega$ is a $\mu$-a.e., Morse cover of $\Omega$. Let 
$\left\langle S_{i}\right\rangle $ be a $\delta$-fine sequence in $\mathcal{S}_{\Omega}$
such that $\left\langle S_{i}\right\rangle $ is $\mu$-exhausting of $\Omega$
and therefore of $X$. We have for any ordering,
\begin{align*}
\left\vert \sum_{i}\Vert f(x_{i})\Vert\mu(S_{i})-\sum_{i}\Vert G(S_{i}
)\Vert\right\vert  &  \leq\sum_{i}\left\vert \ \Vert f(x_{i})\mu(S_{i}
)\Vert-\Vert G(S_{i})\Vert\right\vert \\
&  \leq\sum_{i}\left\Vert f(x_{i})\mu(S_{i})-G(S_{i})\right\Vert
<\frac{\varepsilon}{2}.
\end{align*}
Moreover, for each $j$, 
$\mu\left(  T_{j}\diagdown
{\textstyle\bigcup_{S_{i}\subseteq T_{j}}} S_{i}\right)  =0$, so
\begin{align*}
M  &  \geq\sum_{i}\Vert G(S_{i})\Vert=\sum_{j}\sum_{S_{i}\subseteq T_{j}}\Vert
G(S_{i})\Vert\\
&  \geq\sum_{j}\left\Vert G\left(
{\textstyle\bigcup_{S_{i}\subseteq T_{j}}}
S_{i}\right)  \right\Vert =\sum_{j}\left\Vert G\left(  T_{j}\right)
\right\Vert >M-\frac{\varepsilon}{2}.
\end{align*}
Therefore,
\[
\left\vert \sum_{i}\Vert f(x_{i})\Vert\mu(S_{i})-\ M\right\vert <\varepsilon.
\]
By Theorem~12 in \cite{loebtalvila}, $\Vert f\Vert$ is $\mu$-integrable on
$\Omega$, and therefore on $X$, whence $f$ is $\mu$-integrable on $X$.

To show that $\int_{X}f\,d\mu=G(X)$, we fix a gauge $\delta$ so that Equations
\ref{localerrorsum} and \ref{riemannsum} hold for any sequence 
$\left\langle S_{i}(x_{i})\right\rangle $ of $\delta$-fine sets from 
$\mathcal{S}$ that is
$\mu$-exhausting of $X$. Then given such a sequence, we have
\begin{align*}
&  \left\Vert G(X)-\int_{X}f\,d\mu\right\Vert \\
&  =\left\Vert \sum_{i}\left(  G(S_{i})-\int_{S_{i}}f\,d\mu\right)
\right\Vert \leq\sum_{i}\left\Vert G(S_{i})-\int_{S_{i}}f\,d\mu\right\Vert \\
&  \leq\sum_{i}\left\Vert \int_{S_{i}}f\,d\mu-f(x_{i})\mu(S_{i})\right\Vert
+\sum_{i}\left\Vert f(x_{i})\mu(S_{i})-G(S_{i})\right\Vert <2\varepsilon.
\end{align*}
Since $\varepsilon$ is arbitrary, $\int_{X}f\,d\mu=G(X)$. $\square$

\end{document}